\documentclass[preprint,review,a4paper,11pt,number,sort&compress]{elsarticle}
\usepackage{lscape}
\usepackage{amsmath,amssymb,latexsym}
\usepackage{euscript}
\usepackage[LGR,T1]{fontenc}
\usepackage{bm}

\DeclareSymbolFont{upgreek}{LGR}{cmr}{m}{n}
\SetSymbolFont{upgreek}{bold}{LGR}{cmr}{bx}{n}
\DeclareMathSymbol{\upalpha}{\mathord}{upgreek}{`a}
\DeclareMathSymbol{\upbeta}{\mathord}{upgreek}{`b}
\DeclareMathSymbol{\upgamma}{\mathord}{upgreek}{`g}
\DeclareMathSymbol{\updelta}{\mathord}{upgreek}{`d}
\DeclareMathSymbol{\upepsilon}{\mathord}{upgreek}{`e}
\DeclareMathSymbol{\upzeta}{\mathord}{upgreek}{`z}
\DeclareMathSymbol{\upeta}{\mathord}{upgreek}{`h}
\DeclareMathSymbol{\uptheta}{\mathord}{upgreek}{`j}
\DeclareMathSymbol{\upiota}{\mathord}{upgreek}{`i}
\DeclareMathSymbol{\upkappa}{\mathord}{upgreek}{`k}
\DeclareMathSymbol{\uplambda}{\mathord}{upgreek}{`l}
\DeclareMathSymbol{\upmu}{\mathord}{upgreek}{`m}
\DeclareMathSymbol{\upnu}{\mathord}{upgreek}{`n}
\DeclareMathSymbol{\upxi}{\mathord}{upgreek}{`x}
\DeclareMathSymbol{\upomicron}{\mathord}{upgreek}{`o}
\DeclareMathSymbol{\uppi}{\mathord}{upgreek}{`p}
\DeclareMathSymbol{\uprho}{\mathord}{upgreek}{`r}
\DeclareMathSymbol{\upsigma}{\mathord}{upgreek}{`s}
\DeclareMathSymbol{\uptau}{\mathord}{upgreek}{`t}
\DeclareMathSymbol{\upupsilon}{\mathord}{upgreek}{`u}
\DeclareMathSymbol{\upphi}{\mathord}{upgreek}{`f}
\DeclareMathSymbol{\upchi}{\mathord}{upgreek}{`q}
\DeclareMathSymbol{\uppsi}{\mathord}{upgreek}{`y}
\DeclareMathSymbol{\upomega}{\mathord}{upgreek}{`w}

\DeclareMathSymbol{\Upalpha}{\mathord}{upgreek}{`A}
\DeclareMathSymbol{\Upbeta}{\mathord}{upgreek}{`B}
\DeclareMathSymbol{\Upgamma}{\mathord}{upgreek}{`G}
\DeclareMathSymbol{\Updelta}{\mathord}{upgreek}{`D}
\DeclareMathSymbol{\Upepsilon}{\mathord}{upgreek}{`E}
\DeclareMathSymbol{\Upzeta}{\mathord}{upgreek}{`Z}
\DeclareMathSymbol{\Upeta}{\mathord}{upgreek}{`H}
\DeclareMathSymbol{\Uptheta}{\mathord}{upgreek}{`J}
\DeclareMathSymbol{\Upiota}{\mathord}{upgreek}{`I}
\DeclareMathSymbol{\Upkappa}{\mathord}{upgreek}{`K}
\DeclareMathSymbol{\Uplambda}{\mathord}{upgreek}{`L}
\DeclareMathSymbol{\Upmu}{\mathord}{upgreek}{`M}
\DeclareMathSymbol{\Upnu}{\mathord}{upgreek}{`N}
\DeclareMathSymbol{\Upxi}{\mathord}{upgreek}{`X}
\DeclareMathSymbol{\Upomicron}{\mathord}{upgreek}{`O}
\DeclareMathSymbol{\Uppi}{\mathord}{upgreek}{`P}
\DeclareMathSymbol{\Uprho}{\mathord}{upgreek}{`R}
\DeclareMathSymbol{\Upsigma}{\mathord}{upgreek}{`S}
\DeclareMathSymbol{\Uptau}{\mathord}{upgreek}{`T}
\DeclareMathSymbol{\Upupsilon}{\mathord}{upgreek}{`U}
\DeclareMathSymbol{\Upphi}{\mathord}{upgreek}{`F}
\DeclareMathSymbol{\Upchi}{\mathord}{upgreek}{`Q}
\DeclareMathSymbol{\Uppsi}{\mathord}{upgreek}{`Y}
\DeclareMathSymbol{\Upomega}{\mathord}{upgreek}{`W}

\usepackage{etoolbox}
\makeatletter
\patchcmd{\ps@pprintTitle}
{\let\@oddhead\@empty}
{\def\@oddhead{\mbox{}\hfill\thepage}}
{}{}
\makeatother
\newtheorem{theorem}{Theorem}
\newtheorem{lemma}{Lemma} 
\newtheorem{cor}{Corollary}
\newdefinition{rmk}{Remark}
\newproof{pf}{\bf Proof}
\newproof{pot}{Proof of Theorem \ref{thm2}}
\newcommand{\IR}{{\sf I\! R}}
\newcommand{\ci}{{\mathcal{I}}}
\newcommand{\cj}{{\mathcal{J}}}
\newcommand{\cp}{{\mathcal{P}}}

\begin{document}
\begin{frontmatter}
\title{A simple parallelizable method for the approximate solution of a quadratic transportation problem of large dimension with additional constraints\tnoteref{t1}}
\tnotetext[t1]{This document is a collaborative effort.}
\author[rvt]{S.V.~Rotin}
\ead{s.rotin@gmail.com}
\author[focal]{I.V.~Gusakov\corref{cor1}\fnref{f1}}
\ead{gusakov@goodsforecast.com}
\author[focal]{V.Ya.~Gusakov}
\cortext[cor1]{Corresponding author.}
\fntext[fn1]{The algorithm was proposed and numerically tested by I.V.~Gusakov in 2004.}
\address[rvt]{Priory-Team Information Systems Ltd, Hahoma St., 12, 7565512 Rishon Lezion, Israel}
\address[focal]{GoodsForecast Ltd, Nauchnii Proezd, 19, 117246 Moscow, Russia}
\begin{abstract}
Complexity of the Operations Research Theory tasks can be often diminished in cases that do not require finding the exact solution. For example, forecasting two-dimensional hierarchical time series leads us to the transportation problem with a quadratic objective function and with additional constraints. While solving this task there is no need to minimize objective function with high accuracy, but it is very important to meet all the constraints. In this article we propose a simple iterative algorithm, which can find a valid transportation flow matrix in a limited number of steps while allowing massively parallel computing. 
Method's convergence was studied: a convergence criterion was indicated, as well as the solution's accuracy estimation technique. It was proved that the method converges with the speed of geometric progression, whose ratio weakly depends on the problem's dimension. Numerical experiments were performed to demonstrate the method's efficiency for solving specific large scale transportation problems.
 
\end{abstract}
\end{frontmatter}

\section{Introduction}
It is well known from the Computational Complexity Theory that for some NP-hard optimization problems the complexity class can be reduced in cases when discovering exact global optimum is not necessary \cite{AVL85}. For example, it can be done when it is possible to find a special heuristic algorithm which solves the optimization problem approximately (with an accuracy up to $\varepsilon$ according to the objective function), but with polynomial complexity \cite{Heur2009}. Moreover, for tasks of P complexity class a similar situation also takes place: if it is not necessary to minimize objective function with very high degree of precision, it is often possible to encounter heuristic with a simple structure, allowing as a consequence its widespread application without the necessity to use special mathematical packages. Sometimes this heuristic demonstrates additional benefits, such as high convergence speed as well as parallelizability. However, the convergence conditions analysis even for very simple methods can be far from trivial, while without these conditions application of the method to the real-life problems might be questionable.

In this article we propose to use aforementioned approach to the quadratic transportation problem (QTP) with additional constraints. This particular problem has arisen while working on the task of forecasting two-dimensional hierarchical time-series and was initially posed in \cite{RH2014}. In many businesses data is organized in a hierarchical way, and often has several dimensions. For example, sales transactions between a company and its clients might be represented as a two-dimensional hierarchical data structure (cube), one dimension being product dimension and another client one (note that time dimension is excluded from the cube as it is reserved for the forecasting purpose). Both products and clients are organized in hierarchies: e.g. clients are aggregated by geography while products by brands. For the sake of simplicity, we will further assume that both hierarchies are comprised of only two levels: top level (all products or clients) and bottom level (particular products or clients), but in fact, the method we propose can be repeatedly applied for each level of multi-level hierarchies as well. 

To forecast the data-cube described above, several approaches are applicable. Clearly, all low-level forecasts can be obtained independently. Unfortunately, low-level time series are often quite volatile and forecasting quality can be poor. Another approach is to forecast upper-level time series, and then prorate obtained forecast to the low levels. All pros and cons of these approaches are described in details in \cite{SC2016}. For one-dimensional time series in \cite{RH2011} a special combination approach was proposed, where forecasts are created simultaneously at all levels of the hierarchy, and then optimally reconciled using a regression model. In this article, the reconciliation approach is applied to the two-dimensional time series which leads us to QTP as described below.

Let's consider a structured set of random functions and their partial realizations, which represent actual sales of a company's products to its clients (the first index indicates client and the second one -- product):

\begin{equation}
A_{ij}(t)\, ,\;\; t_0 < t\quad\mbox{and}\quad
a_{ij}(t)\, ,\;\; t_0 < t \le t_e\, ,\qquad
i=\overline{1,\,n}\, ,\;\; j=\overline{1,\,m}\, ,
\label{Aa}
\end{equation}

Also, additional random functions, strictly dependent on original ones, are taken into consideration:
\begin{eqnarray}
\mbox{}\qquad & \displaystyle S(t) = \sum\limits_{i=1}^{n} \sum\limits_{j=1}^{m} A_{ij}(t)\, ,\hphantom{-------}
& \label{S}\\
\mbox{}\qquad & \displaystyle P_i (t) = \sum\limits_{j=1}^{m} A_{ij}(t)\, ,\;\; i=\overline{1,\,n}\, ,\qquad
Q_j (t) = \sum\limits_{i=1}^{n} A_{ij}(t)\, ,\;\; j=\overline{1,\,m}\, .& \label{PQ}
\end{eqnarray}
Suppose that at a particular time $t>t_e$ we have obtained independent forecasts: $\{a_{ij}(t)\}$, $\{p_{i}(t)\}$, $\{q_{j}(t)\}$ and $s(t)$. Further on variable $t$ in all formulas will be omitted for the sake of simplicity. 

In the first step, we need to balance upper-level forecasts, i.e. to make sure that:
$s = \sum\limits_{i=1}^n p_i = \sum\limits_{j=1}^m q_j$.
In many businesses, this condition is satisfied automatically due to the nature of the top-down forecasting process: forecasts $p_i$ and $q_j$ are in fact not independent but obtained as a shares of forecast $s$. In case this condition is not satisfied, a relatively simple optimization task might be solved to perform optimal balancing. This optimization task will be described in details in our article to follow. 

Moreover, it is clear that upper-level forecast cannot differ dramatically from the sum of lower-level forecasts, so we assume that
$s \approx \sum\limits_{j=1}^m\sum\limits_{i=1}^n a_{ij}$. We now need to find corrected sales forecast matrix $\bm X$, balanced by both rows and columns, while at the same time minimally different from the original matrix $\bm A=\left(a_{ij}\right)$. We assume that both actual and future sales are non-negative (negative sales are nothing but returns of the products, which are not taken into account when performing forecasting task). Moreover, often there are products that are not sold to particular clients due to a company policy or other factors, which means that if there are zeros in some positions of the matrix $\bm A$, these zeros should also be present in the same positions of the matrix $\bm X$. As a consequence, we obtain the following quadratic programming problem:
\begin{eqnarray}
& \displaystyle
J^2(\bm X) \equiv\|\bm X - \bm A\|^2_2 = \sum\limits_{i\in\ci}\sum\limits_{j\in\cj} \left(x_{ij} - a_{ij}\right)^2
\rightarrow\min_{\bm X \in \IR^{n\times m}} &
\label{J}\\[5.pt]
& \displaystyle\!\!\!\!\!\!\!\!\!\!\!\!\!\! 
\sum\limits_{j\in\cj} x_{ij} = p_i\, ,\quad i\in\ci\, , &
\label{Xp}\\
& \displaystyle\!\!\!\!\!\!\!\!\!\!\!\!\!
\sum\limits_{i\in\ci} x_{ij} = q_j\, ,\quad j\in\cj\, , &
\label{Xq}\\
& \displaystyle\!\!\!\!\! 
x_{ij} \ge 0\, ,\quad (i,j)\in\ci\times\cj\, , &
\label{Neotric}\\
& \displaystyle \quad\,
x_{ij} \le 0\, ,\quad (i,j)\in\cp\subset\ci\times\cj\, . &
\label{Dop}
\end{eqnarray}
with the following sets of indexes:
\begin{equation}
\ci=\{1,\,2,...,\,n\},\qquad \cj=\{1,\,2,...,\,m\} ,\qquad
\cp =\{(i,j)\in\ci\times\cj: a_{ij} = 0\}\, .
\label{IJP}
\end{equation}

Problem (\ref{J}-\ref{Dop}) is quadratic transportation problem (QTP) with additional constraints (\ref{Dop}). Without going into details of the verification procedure \cite{Polyak}, we will further assume that the full system of restrictions is consistent.

\section{Numerical algorithm and its convergence analysis} 
\noindent
Suppose $\bm X^{(0)} = \bm A$. For $s=1, 3, 5, ...$ let's consider the following two-step iterative procedure:
\begin{eqnarray}
&& x^{(s)}_{ij} = \frac{p_i\, x^{(s-1)}_{ij}}{\sum_{l=1}^{m}x^{(s-1)}_{il}} \, ,\qquad (i,j)\in\ci\times\cj \, , \label{s}\\
&& x^{(s+1)}_{ij} = \frac{q_j\, x^{(s)}_{ij}}{\sum_{k=1}^{n}x^{(s)}_{kj}} \, ,\qquad (i,j)\in\ci\times\cj \, . \label{sp1}
\end{eqnarray}
We will further use term "iteration"\ to designate elementary periodical unit of the aforementioned iterative procedure and term "step"\ to designate first and second steps inside every iteration. Nevertheless, we will use continuous numbering of the steps through all the iterations, so $s=1, 2$ refer to the steps of the first iteration, $s=3, 4$ -- to the steps of the second iteration and so on.

At every step of a particular iteration balancing of the forecast matrix either by lines or by columns is performed. It's quite clear that restrictions (\ref{Neotric}-\ref{Dop}) in this case are fulfilled automatically. As we start our iterative procedure from matrix $\bm A$, there is hope that solutions, obtained by this heuristic, will provide the objective function (\ref{J}) with volumes, relatively close to the minimal ones. Mind, that as we could have used other norm (not Euclidean) to define the measure of the difference between the objective matrix and the original one, there is no need to minimize the Euclidean norm exactly. What is indeed important, is for the solution $\bm X$ to be in some way close to the original matrix $\bm A$.

\begin{theorem}
	\ \ Matrix sequence $\left\lbrace \bm X^{(s)}\right\rbrace $ converges to the admissible set, defined by the restrictions  (\ref{Xp}-\ref{Xq}), with the speed of geometric progression, if the original data satisfy the following conditions:
	\begin{itemize}
		\item[a)] matrix $\bm A\in\IR_{+}^{n\times m}$ does not contain neither zero lines no zero columns,
		\item[b)] all the elements of vectors $\bm p\in\IR_{+}^{n}$ and $\bm q\in\IR_{+}^{m}$ are strictly positive,
		\item[c)] $\epsilon_p(\bm A)>0$ and for the initial approximation $\bm X^{(0)}$ the following is true: 
		\begin{equation}
		V(\bm X^{(0)})\left(1+O\!\left(\frac{V(\bm X^{(0)})}{\min_{i\in\ci} p_i}\right)\right) \le \frac{1}{2}\,\epsilon_p(\bm X^{(0)})\min_{i\in\ci} p_i\, , \label{OblSx}
		\end{equation}		
		where
		\begin{equation}
		\epsilon_p(\bm X) \equiv 
		\min_{\tilde{\ci}\subset\ci} \min_{\tilde{\cj}\subset\cj}\left\{
		\min_{j\in\tilde{\cj}}\sum\limits_{i\in\ci\setminus\tilde{\ci}}\frac{x_{ij}}{\sum\limits_{k=1}^n x_{kj}} + 
		\min_{j\in\cj\setminus\tilde{\cj}}\sum\limits_{i\in\tilde{\ci}}\frac{x_{ij}}{\sum\limits_{k=1}^n x_{kj}}\right\}\, ,
		\label{eps}
		\end{equation}
		\begin{equation}
		V(\bm X) \equiv 
		\sum\limits_{i=1}^n\left|\sum\limits_{j=1}^m x_{ij} - p_i\,\right| +
		\sum\limits_{j=1}^m\left|\sum\limits_{i=1}^n x_{ij} - q_j\,\right| \ge 0 \, .
		\label{Lyapunov}
		\end{equation}
	\end{itemize}
\end{theorem}

\begin{pf} Let's take (\ref{Lyapunov}) as a Lyapunov function \cite{Polyak}. For any matrix $\tilde{\bm X}$, which belongs to the admissible set (\ref{Xp}-\ref{Xq}), and only for such matrices, $V(\tilde{\bm X})=0$. Let's show that the results of iterations of the method (\ref{s}-\ref{sp1}) does not move away from the admissible set, that is $V(\bm X^{(s+1)}) \le V(\bm X^{(s)}) \le V(\bm X^{(s-1)})$. For the second step (\ref{sp1}) we have the following upper estimate:
\begin{eqnarray}
& V(\bm X^{(s+1)}) = \sum\limits_{i=1}^n \left|\sum\limits_{j=1}^{m}x^{(s+1)}_{ij} - p_i\,\right| + 0 = 
\sum\limits_{i=1}^n\left|\sum\limits_{j=1}^m \frac{q_j\, x^{(s)}_{ij}}{\sum_{k=1}^{n}x^{(s)}_{kj}} - p_i\,\right| = {}
& \nonumber\\
& {} = \sum\limits_{i=1}^n\left|\sum\limits_{j=1}^m \frac{q_j\, x^{(s)}_{ij}}{\sum_{k=1}^{n}x^{(s)}_{kj}} - \sum\limits_{j=1}^m x^{(s)}_{ij} + \sum\limits_{j=1}^m x^{(s)}_{ij} - p_i\,\right| \le
\sum\limits_{i=1}^n\left(\left|\sum\limits_{j=1}^m \frac{q_j\, x^{(s)}_{ij}}{\sum_{k=1}^{n}x^{(s)}_{kj}} - \sum\limits_{j=1}^m x^{(s)}_{ij}\right| + \left|\sum\limits_{j=1}^m x^{(s)}_{ij} - p_i\,\right|\right)
= {} & \nonumber\\
& {} =\sum\limits_{i=1}^n\left(\left|\sum\limits_{j=1}^m x^{(s)}_{ij} - p_i\,\right| + \left|\sum\limits_{j=1}^m \left(q_j - \sum\limits_{k=1}^{n}x^{(s)}_{kj}\right) \frac{ x^{(s)}_{ij}}{\sum_{k=1}^{n}x^{(s)}_{kj}}\right| \right) \le
\sum\limits_{i=1}^n \left|\sum\limits_{j=1}^{m}x^{(s)}_{ij} - p_i\,\right| +  {}& \nonumber \\
&{} + \sum\limits_{i=1}^n\sum\limits_{j=1}^m \!\left|\, q_j - \sum\limits_{k=1}^{n}x^{(s)}_{kj}\right| \left|\frac{ x^{(s)}_{ij}}{\sum_{k=1}^{n}x^{(s)}_{kj}}\right| =
\sum\limits_{i=1}^n \left|\sum\limits_{j=1}^{m}x^{(s)}_{ij} - p_i\,\right| +
\sum\limits_{j=1}^m \left|\sum\limits_{i=1}^{n}x^{(s)}_{ij} - q_j\right| \sum\limits_{i=1}^n \frac{ \left|\, x^{(s)}_{ij}\right|}{\left|\sum_{k=1}^{n}x^{(s)}_{kj}\right|}
= V(\bm X^{(s)})& \nonumber
\end{eqnarray}
For the first step (\ref{s}) estimation is made in a similar way. The last equality in the previous formula is the consequence of the fact that all elements of matrices $\bm X^{(s)},\; s=1,2,3,...$ are non-negative. The last inequality in this formula might be strict. To understand when exactly it happens, let's separately consider the expression below  
\begin{eqnarray}
& L(\bm X^{(s)}) \equiv \sum\limits_{i=1}^n\left|\sum\limits_{j=1}^m \left(q_j - \sum\limits_{k=1}^{n}x^{(s)}_{kj}\right) \frac{ x^{(s)}_{ij}}{\sum_{k=1}^{n}x^{(s)}_{kj}}\right| = \sum\limits_{i\in\ci_{+}}\!\left(\sum\limits_{j\in\cj_{+}} \left|\sum\limits_{k=1}^{n}x^{(s)}_{kj}\!-q_j\right| \frac{ x^{(s)}_{ij}}{\sum_{k=1}^{n}x^{(s)}_{kj}} - {}\right.& \nonumber\\[10pt]
&\left.{}-\sum\limits_{j\in\cj_{-}} \left|\sum\limits_{k=1}^{n}x^{(s)}_{kj}\!-q_j\right| \frac{ x^{(s)}_{ij}}{\sum_{k=1}^{n}x^{(s)}_{kj}}\right) + 
\sum\limits_{i\in\ci_{-}}\!\left(\sum\limits_{j\in\cj_{-}} \left|\sum\limits_{k=1}^{n}x^{(s)}_{kj}\!-q_j\right| \frac{ x^{(s)}_{ij}}{\sum_{k=1}^{n}x^{(s)}_{kj}} - \sum\limits_{j\in\cj_{+}} \left|\sum\limits_{k=1}^{n}x^{(s)}_{kj}\!-q_j\right| \frac{ x^{(s)}_{ij}}{\sum_{k=1}^{n}x^{(s)}_{kj}}\right). & \nonumber
\end{eqnarray}
Here we use the fact that sum of the elements of sequence $\left\{q_j - \sum_{k=1}^{n}x^{(s)}_{kj}\right\}$ by all $j$ at every iteration is equal to zero, and, as a consequence, the set of indexes is divided into two disjoint subsets: $\cj=\cj_{+}\cup\cj_{-}$, that is the subset of indexes where elements of the sequence are positive, and the subset of indexes where they are negative (both subsets are not empty). But in this case, the set of indexes of the external sum is also divided into two subsets: $\ci=\ci_{+}\cup\ci_{-}$, where the first subset corresponds with positive internal sums inside the module sign, while the second subset -- with negative ones. Let's show that both these subsets cannot be empty as well:
\begin{equation}
\sum\limits_{i=1}^n\left(\sum\limits_{j=1}^m \left(q_j - \sum\limits_{k=1}^{n}x^{(s)}_{kj}\right) \frac{ x^{(s)}_{ij}}{\sum_{k=1}^{n}x^{(s)}_{kj}}\right) = \sum\limits_{j=1}^m \left(q_j - \sum\limits_{k=1}^{n}x^{(s)}_{kj}\right) \sum\limits_{i=1}^n\frac{x^{(s)}_{ij}}{\sum_{k=1}^{n}x^{(s)}_{kj}} = \sum\limits_{j=1}^m q_j - \sum\limits_{j=1}^m\sum\limits_{i=1}^{n}x^{(s)}_{ij} = 0
\, . \nonumber
\end{equation}
Let's now in the expression for $L(\bm X^{(s)})$ add and subtract those 2 terms out of 4, which participate in this expression with preceding minus sign:
\begin{eqnarray}
& L(\bm X^{(s)}) = \sum\limits_{i=1}^n\sum\limits_{j=1}^m \left|\sum\limits_{k=1}^{n}x^{(s)}_{kj}\!-q_j\right| \frac{ x^{(s)}_{ij}}{\sum_{k=1}^{n}x^{(s)}_{kj}} - 2 \sum\limits_{i\in\ci_{+}}\!\left(\sum\limits_{j\in\cj_{-}} \left|\sum\limits_{k=1}^{n}x^{(s)}_{kj}\!-q_j\right| \frac{ x^{(s)}_{ij}}{\sum_{k=1}^{n}x^{(s)}_{kj}}\right) - {} &\nonumber\\[10pt]
& {} - 2 \sum\limits_{i\in\ci_{-}}\!\left(\sum\limits_{j\in\cj_{+}} \left|\sum\limits_{k=1}^{n}x^{(s)}_{kj}\!-q_j\right| \frac{ x^{(s)}_{ij}}{\sum_{k=1}^{n}x^{(s)}_{kj}}\right) =
\sum\limits_{j=1}^m \left|\sum\limits_{k=1}^{n}x^{(s)}_{kj}\!-q_j\right| \sum\limits_{i=1}^n \frac{ x^{(s)}_{ij}}{\sum_{k=1}^{n}x^{(s)}_{kj}} -
2 \sum\limits_{j\in\cj_{-}} \left|\sum\limits_{k=1}^{n}x^{(s)}_{kj}\!-q_j\right|\times {} & \nonumber\\[10pt]
& {}\times\sum\limits_{i\in\ci_{+}} \frac{ x^{(s)}_{ij}}{\sum_{k=1}^{n}x^{(s)}_{kj}} -
2 \sum\limits_{j\in\cj_{+}} \left|\sum\limits_{k=1}^{n}x^{(s)}_{kj}\!-q_j\right|\times\sum\limits_{i\in\ci_{-}} \frac{ x^{(s)}_{ij}}{\sum_{k=1}^{n}x^{(s)}_{kj}} \le \sum\limits_{j=1}^m \left|\sum\limits_{k=1}^{n}x^{(s)}_{kj}\!-q_j\right|  - 2 \left(\sum\limits_{j\in\cj_{-}} \left|\sum\limits_{k=1}^{n}x^{(s)}_{kj}\!-q_j\right|\times {} \right. & \nonumber\\[10pt]
& \left. {}\times\min\limits_{j\in\cj_{-}}\sum\limits_{i\in\ci_{+}} \frac{ x^{(s)}_{ij}}{\sum_{k=1}^{n}x^{(s)}_{kj}} +
\sum\limits_{j\in\cj_{+}} \left|\sum\limits_{k=1}^{n}x^{(s)}_{kj}\!-q_j\right|\times\min\limits_{j\in\cj_{+}}\sum\limits_{i\in\ci_{-}} \frac{ x^{(s)}_{ij}}{\sum_{k=1}^{n}x^{(s)}_{kj}}\right) \le \sum\limits_{j=1}^m \left|\sum\limits_{k=1}^{n}x^{(s)}_{kj}\!-q_j\right| - {}& \nonumber\\[10pt]
& \sum\limits_{j=1}^m \left|\sum\limits_{k=1}^{n}x^{(s)}_{kj}\!-q_j\right|\!\left(\min\limits_{j\in\cj_{-}}\sum\limits_{i\in\ci_{+}} \frac{ x^{(s)}_{ij}}{\sum_{k=1}^{n}x^{(s)}_{kj}} +
\min\limits_{j\in\cj_{+}} \sum\limits_{i\in\ci_{-}} \frac{ x^{(s)}_{ij}}{\sum_{k=1}^{n}x^{(s)}_{kj}} \right) \equiv R(\bm X^{(s)}) .& \nonumber
\end{eqnarray}
The last inequality is correct because for 4 non-negative numbers the following implication is true:
\begin{equation}
\upsilon_{+} = \upsilon_{-} \quad\Rightarrow\quad
2(\upsilon_{-}\mu_{-} + \upsilon_{+}\mu_{+}) = 2\upsilon_{\pm}(\mu_{-} + \mu_{+}) = (\upsilon_{-} + \upsilon_{+})(\mu_{-} + \mu_{+}) \, . 
\nonumber\label{Schema}
\end{equation}
Let's make the first estimation of the proof stronger, taking into account that 
$\sum\limits_{i=1}^n \left|\sum\limits_{j=1}^{m}x^{(s)}_{ij} - p_i\,\right| = 0$:
\begin{equation}
V(\bm X^{(s+1)}) \le 0 + L(\bm X^{(s)}) \le R(\bm X^{(s)})
\le \left(1 - \epsilon_p(\bm X^{(s)})\right) V(\bm X^{(s)})\, . \label{Veval}
\end{equation}
The last inequality is necessary because we need to obtain iteration-independent estimation. As current division of the sets of indexes into subsets for the iteration is not known, we have to use double minimum by all possible subsets  $\tilde{\ci}\subset\ci$, $\tilde{\cj}\subset\cj$. This continues the upper estimation.

Suppose now, that the initial approximation  $\bm X^{(0)}$ is not in the "dead zone"\ ($\epsilon_p(\bm X^{(0)})\ne 0$) and is close enough to the admissible set for the condition (\ref{OblSx}) to be true. For the odd $s$ and arbitrary subsets $\tilde{\ci}$ and $\tilde{\cj}$ we can estimate the following expression:
\begin{eqnarray}
&\displaystyle G(\tilde{\ci},\tilde{\cj},\bm X^{(s)}) \equiv \min\limits_{j\in\tilde{\cj}} \frac{\sum\limits_{i\in\tilde{\ci}} x^{(s)}_{ij}}{\sum\limits_{k=1}^{n}x^{(s)}_{kj}} =
\min\limits_{j\in\tilde{\cj}} \frac{\sum\limits_{i\in\tilde{\ci}} \frac{p_i\, x^{(s-1)}_{ij}}{\sum_{l=1}^{m}x^{(s-1)}_{il}-p_i+p_i}}{\sum\limits_{i=1}^{n}\frac{p_i\, x^{(s-1)}_{ij}}{\sum_{l=1}^{m}x^{(s-1)}_{il}-p_i+p_i}} \ge 
\min\limits_{j\in\tilde{\cj}} \frac{\min\limits_{i\in\tilde{\ci}} \frac{p_i}{\sum_{l=1}^{m}x^{(s-1)}_{il}-p_i+p_i}}{\max\limits_{i\in\ci}\frac{p_i}{\sum_{l=1}^{m}x^{(s-1)}_{il}-p_i+p_i}}\, \frac{\sum_{i\in\tilde{\ci}} x^{(s-1)}_{ij}}{\sum_{i=1}^{n} x^{(s-1)}_{ij}} \ge {}& \nonumber\\
&\displaystyle{}\ge\frac{\frac{1}{1+\max_{i\in\ci}\left|\sum_{l=1}^{m}x^{(s-1)}_{il}-p_i\right|/p_i}}{\frac{1}{1-\max_{i\in\ci}\left|\sum_{l=1}^{m}x^{(s-1)}_{il}-p_i\right|/p_i}} \min\limits_{j\in\tilde{\cj}} \frac{\sum\limits_{i\in\tilde{\ci}} x^{(s-1)}_{ij}}{\sum\limits_{k=1}^{n} x^{(s-1)}_{kj}} =
\frac{1-\max_{i\in\ci}p_i^{-1} \max_{i\in\ci}\left|\sum_{l=1}^{m}x^{(s-1)}_{il}-p_i\right|
}{1+\max_{i\in\ci}p_i^{-1} \max_{i\in\ci}\left|\sum_{l=1}^{m}x^{(s-1)}_{il}-p_i\right|} \, G(\tilde{\ci},\tilde{\cj},\bm X^{(s-1)}) \ge {}& \nonumber\\
&\displaystyle{}\ge\frac{\min_{i\in\ci} p_i - \sum_{i=1}^{n}\left|\sum_{l=1}^{m}x^{(s-1)}_{il}-p_i\right|}{\min_{i\in\ci} p_i + \sum_{i=1}^{n}\left|\sum_{l=1}^{m}x^{(s-1)}_{il}-p_i\right|}\, G(\tilde{\ci},\tilde{\cj},\bm X^{(s-1)}) = {}
\frac{\min_{i\in\ci} p_i - V(\bm X^{(s-1)})}{\min_{i\in\ci} p_i + V(\bm X^{(s-1)})}\, G(\tilde{\ci},\tilde{\cj},\bm X^{(s-1)})\, .&
\nonumber
\end{eqnarray}
It's easy to check that for odd number $s>1$:
\begin{equation}
G(\tilde{\ci},\tilde{\cj},\bm X^{(s-1)}) = G(\tilde{\ci},\tilde{\cj},\bm X^{(s-2)}) \, . \label{ravno}
\end{equation}
If we continue the obtained estimation by all iterations, we get:
\begin{equation}
G(\tilde{\ci},\tilde{\cj},\bm X^{(s)}) \ge
\frac{\min_{i\in\ci} p_i - V(\bm X^{(s-1)})}{\min_{i\in\ci} p_i + V(\bm X^{(s-1)})}\times ... \times\frac{\min_{i\in\ci} p_i - V(\bm X^{(0)})}{\min_{i\in\ci} p_i + V(\bm X^{(0)})}\, G(\tilde{\ci},\tilde{\cj},\bm X^{(0)})\, .\label{Times}
\end{equation}
In total we have $(s+1)/2$ multipliers. Now we substitute (\ref{Times}) into (\ref{Veval}) with the relevant subsets of indexes, and taking into account that Lyapunov function is non-increasing, we receive:
\begin{eqnarray}
& V(\bm X^{(s+1)})  & \le
\left(1 - \min_{\tilde{\ci}\subset\ci} \min_{\tilde{\cj}\subset\cj}
\left\{G(\ci\setminus\tilde{\ci},\,\tilde{\cj},\,\bm X^{(s)}) + 
G(\tilde{\ci},\,\cj\setminus\tilde{\cj},\,\bm X^{(s)})\right\}\right) V(\bm X^{(s)}) \le {} \nonumber\\
&& {} \le \left(1 - \epsilon_p(\bm X^{(0)})\!\!\prod\limits_{t=0}^{(s-1)/2}\frac{\min_{i\in\ci} p_i - V(\bm X^{(2t)})}{\min_{i\in\ci} p_i + V(\bm X^{(2t)})}\right) V(\bm X^{(s-1)})\, . \label{Veval2}
\end{eqnarray}
In accordance with Lemma~\ref{Lem1}, sequence $V(\bm X^{(0)}),\,V(\bm X^{(2)}),\,...,V(\bm X^{(s+1)})$, which satisfy (\ref{Veval2}), converges to zero with the speed of geometric progression, if the initial approximation satisfy the condition (\ref{OblSx}), which corresponds with inequality (\ref{omal}) in the lemma. $\Box$
\end{pf}

\begin{lemma}\label{Lem1}
	\ \ Let a sequence of positive numbers $\{V_i\}_{i=0}^{\infty}$ be given by the recurrence relation
		\begin{equation}
	V_i = \left(1 - \varepsilon_1 \prod\limits_{j=0}^{i-1} 
	\frac{1 - \varepsilon_2 V_j}{1 + \varepsilon_2 V_j}\right) V_{i-1}\, , \qquad i=1,\, 2,\, 3,\, ... \label{Rek}
	\end{equation}
	then there does exist a majorizing geometric progression $\{U_i\}_{i=0}^{\infty}$ with 	
	ratio 
	$(1-\alpha)$ (between zero and one) so that: 
	\begin{equation}
	V_0 = U_0\, , \qquad 0 < V_i \le U_i\, , \qquad i=1,\, 2,\, 3,\, ... \label{Mazhor}
	\end{equation} 
if for the initial term of the original sequence following conditions are satisfied:
	\begin{eqnarray}
	& 0 \le \varepsilon_2 V_0 < 1\, , & \label{e2}\\[10pt]
	& \displaystyle 0 < \varepsilon_1 < \frac{1+\varepsilon_2 V_0}{1-\varepsilon_2 V_0}\, ,   & \label{e1}\\[10pt]
	& 2 \varepsilon_2 V_0\left(\vphantom{2^{2^2}}1+O\left(\varepsilon_2 V_0\right)\right) \le \varepsilon_1\, .& \label{omal}
	\end{eqnarray}
\end{lemma}

\begin{pf} Let's rewrite recurrence relation (\ref{Rek}) in the following way:
\begin{equation}
\ln\left(-\frac{\Delta V_{i}}{V_{i}}\right) = \ln\varepsilon_1 + \sum\limits_{j=0}^{i} \ln (1 - \varepsilon_2 V_j) - \sum\limits_{j=0}^{i} \ln (1 + \varepsilon_2 V_j)\, , \qquad i=0,\, 1,\, 2,\, ... \label{Rek2}
\end{equation}
where $\Delta W_i \equiv W_{i+1} - W_i$ is \textit{1st order finite difference operator} on an arbitrary numeric sequence $\{W_i\}_{i=0}^{\infty}$. It can be shown that this operator has the following property:
\begin{equation}
\Delta\ln(W_i) = \frac{\Delta W_i}{W_{i+1}} + \omega\left(-\frac{\Delta W_i}{W_{i+1}}\right)\, , \qquad i=0,\, 1,\, 2,\, ... \label{Difln}
\end{equation}
where $\omega(w)\equiv w-\ln(1+w)$. Let's now apply finite difference operator to the equation (\ref{Rek2}):
\begin{equation}
\frac{\Delta\left(-\frac{\Delta V_{i-1}}{V_{i-1}}\right)}{-\frac{\Delta V_{i}}{V_{i}}} + \omega\!\left(-\frac{\Delta\left(-\frac{\Delta V_{i-1}}{V_{i-1}}\right)}{-\frac{\Delta V_{i}}{V_{i}}}\right) = \ln (1 - \varepsilon_2 V_i) - \ln (1 + \varepsilon_2 V_i)\, , \qquad i=1,\, 2,\, 3,\, ... \label{DifRek}
\end{equation}
Under condition (\ref{e1}) this sequence is \textit{strictly
monotone} decreasing and thus (\ref{DifRek}) can be rewritten as:
\begin{equation}
\Delta\left(-\frac{\Delta V_{i-1}}{V_{i-1}}\right) = \frac{\Delta V_{i}}{V_{i}}\,\omega\!\left(-\frac{\Delta\left(-\frac{\Delta V_{i-1}}{V_{i-1}}\right)}{-\frac{\Delta V_{i}}{V_{i}}}\right) + \varepsilon_2 \Delta V_{i}\,\frac{\ln (1 + \varepsilon_2 V_i) - \ln (1 - \varepsilon_2 V_i)}{\varepsilon_2 V_i}\, , \qquad i=1,\, 2,\, 3,\, ... \label{Rek3}
\end{equation}
After summing in (\ref{Rek3}) from 1 to  $i$ and reducing on the left side the same terms we get:
\begin{equation}
\left(-\frac{\Delta V_{i}}{V_{i}} - -\frac{\Delta V_{0}}{V_{0}}\right) = \sum\limits_{j=1}^{i}\left(\frac{\Delta V_{j}}{V_{j}}\,\omega\!\left(-\frac{\Delta\left(-\frac{\Delta V_{j-1}}{V_{j-1}}\right)}{-\frac{\Delta V_{j}}{V_{j}}}\right) + \varepsilon_2 \Delta V_{j}\,\frac{\ln (1 + \varepsilon_2 V_j) - \ln (1 - \varepsilon_2 V_j)}{\varepsilon_2 V_j}\right)\, . \label{Rek4}
\end{equation}
In accordance with Lemma~\ref{Lem2} majorizing geometric
progression $\{U_i\}_{i=0}^{\infty}$ will exist, if
\begin{equation}
-\frac{\Delta V_{0}}{V_{0}} + \sum\limits_{j=1}^{i}\left(\frac{\Delta V_{j}}{V_{j}}\,\omega\!\left(-\frac{\Delta\left(-\frac{\Delta V_{j-1}}{V_{j-1}}\right)}{-\frac{\Delta V_{j}}{V_{j}}}\right) + \varepsilon_2 \Delta V_{j}\,\frac{\ln (1 + \varepsilon_2 V_j) - \ln (1 - \varepsilon_2 V_j)}{\varepsilon_2 V_j}\right) \ge\alpha > 0\, ,\label{Krit}
\end{equation}
for all $i \ge 1$. In this case all the minorants $\psi_i$ are equal to constant $\alpha$ and inequalities (\ref{psi}) are fulfilled, the last one as a consequence of (\ref{e1}). Let's designate expression in the left part of (\ref{Krit}) as $R_i$ and calculate for it an estimate from below:
\begin{eqnarray}
& \displaystyle R_i =  \varepsilon_1\frac{1-\varepsilon_2 V_0}{1+\varepsilon_2 V_0} + \sum\limits_{j=1}^{i}\left(\frac{\Delta V_{j}}{V_{j}}\,\omega\!\left(\frac{1+\varepsilon_2 V_j}{1-\varepsilon_2 V_j} - 1\right) + \varepsilon_2  \Delta V_{j}\,\frac{\ln (1 + \varepsilon_2 V_j) - \ln (1 - \varepsilon_2 V_j)}{\varepsilon_2 V_j}\right) \ge {}&
\nonumber\\[10pt]
& \displaystyle {} \ge\varepsilon_1\frac{1-\varepsilon_2 V_0}{1+\varepsilon_2 V_0} + \sum\limits_{j=1}^{i}\Delta V_{j} \max_{j=1,...,i}\left(\frac{1}{V_{j}}\,\omega\!\left(\frac{2\varepsilon_2 V_j}{1-\varepsilon_2 V_j}\right) + \varepsilon_2\,\frac{\ln (1 + \varepsilon_2 V_j) - \ln (1 - \varepsilon_2 V_j)}{\varepsilon_2 V_j}\right) \ge {} & \nonumber\\[10pt]
& \displaystyle {} \ge\varepsilon_1\frac{1-\varepsilon_2 V_0}{1+\varepsilon_2 V_0} - V_1\left(\frac{1}{V_{1}}\,\omega\!\left(\frac{2\varepsilon_2 V_1}{1-\varepsilon_2 V_1}\right) + \varepsilon_2\,\frac{\ln (1 + \varepsilon_2 V_1) - \ln (1 - \varepsilon_2 V_1)}{\varepsilon_2 V_1}\right) \ge\alpha > 0& \nonumber
\end{eqnarray}
The last inequality is fulfilled only in case the initial term was selected not too far from zero:
\begin{equation}
\left( \frac{2\varepsilon_2 V_0}{1-\varepsilon_2 V_1} \left(1-\varepsilon_1\,\frac{1-\varepsilon_2 V_0}{1+\varepsilon_2 V_0}\right) + \alpha\right)\frac{1+\varepsilon_2 V_0}{1-\varepsilon_2 V_0} \le\varepsilon_1 \, , \label{Zona}
\end{equation}
Performed above simplification is based on (\ref{e2}) and on definition of the function $\omega$. Let the ratio of the geometric progression be defined by the following parameter 
\begin{equation}
\alpha \equiv \varepsilon_1\,\frac{1-\varepsilon_2 V_0}{1+\varepsilon_2 V_0}
\,\frac{2\varepsilon_2 V_0}{1-\varepsilon_2 V_1} = 2\,\varepsilon_1\varepsilon_2 V_0 \left(\vphantom{2^{2^2}} 1+(1+\varepsilon_1)\,\varepsilon_2 V_0\right)\!
\left(\frac{1-\varepsilon_2 V_0}{1+\varepsilon_2 V_0}\right)^{\!\!2} > 0
\, , \label{Alpha}
\end{equation}
then, according to (\ref{Zona}), for the initial term to be close to zero is necessary:
\begin{equation}
\frac{2\varepsilon_2 V_0}{1+(1+\varepsilon_1)\,\varepsilon_2 V_0} 
\left(\frac{1+\varepsilon_2 V_0}{1-\varepsilon_2 V_0}\right)^{\!\!2} =
2\,\varepsilon_2 V_0 \left(\vphantom{2^{2^2}} 1 + \varepsilon_2 V_0\, r(\varepsilon_2 V_0,\varepsilon_1)\right) \le\varepsilon_1 \, , \label{Zona2}
\end{equation}
where auxiliary function $r$ has, taking into account (\ref{e2}), the following form and properties:
\begin{equation}
r(w,\varepsilon_1) \equiv \frac {4-\left (1+\varepsilon_1\right )\left (1-w\right )^{2}}{\left (1+\varepsilon_1\,w-\left (1+\varepsilon_1\right )w^{2}\right )\left (1-w\right )} > 0\, ,\qquad
\lim\limits_{w\to+0} r(w,\varepsilon_1) = 3 - \varepsilon_1
\, . \label{RF}
\end{equation}
As a consequence, from (\ref{Zona2}-\ref{RF}) we get (\ref{omal}). $\Box$ 
\end{pf}

\begin{lemma} \label{Lem2}
	\ \ Let the sequence of functions $\{\varphi_i(v_0,...,v_i)\}_{i=0}^{\infty}$ be positive and does not exceed 1 in the respective domains
	\begin{equation}
	D_i = \{(v_0,...,v_i): 0\le v_i\le\cdots\le v_0\}\, , \qquad i=0,\, 1,\, 2,\, ... \label{DefBer}
	\end{equation} 
	In case there does exist minorizing sequence of functions $\{\psi_i(v_0,...,v_i)\}_{i=0}^{\infty}$  so that
	\begin{equation}
	0 < \psi_i(v_0,...,v_i) \le\varphi_i(v_0,...,v_i)\le 1\ \quad\mbox{in } D_i\, , \qquad i=0,\, 1,\, 2,\, ... \label{psi}
	\end{equation}
	monotone decreasing by the set of its arguments, that is
	\begin{equation}
	v_0 \le u_0 \wedge\cdots\wedge v_i \le u_i\quad\Rightarrow\quad
	\psi_i(v_0,...,v_i) \ge\psi_i(u_0,...,u_i)\, , \qquad i=0,\, 1,\, 2,\, ... \label{Monoton}
	\end{equation}
	then for recurrently defined positive numeric sequence $\{V_i\}_{i=0}^{\infty}$
	\begin{equation}
	\Delta V_i = -\varphi_i(V_0,...,V_i)\, V_i \, , \qquad i=0,\, 1,\, 2,\, ... \label{Vi}
	\end{equation}
    there does exist a majorizing sequence $\{U_i\}_{i=0}^{\infty}$:
	\begin{equation}
	V_0 = U_0\, , \qquad V_i \le U_i \, , \qquad i=1,\, 2,\, 3,\, ... \label{Lmajor}
	\end{equation}
	given by the recurrence relation:
	\begin{equation}
	\Delta U_i = -\psi_i(U_0,...,U_i)\, U_i \, , \qquad i=0,\, 1,\, 2,\, ... \label{Ui}
	\end{equation}
\end{lemma}

\begin{pf} 
Let's apply the method of mathematical induction. For $i=1$:
\begin{equation}
V_1 = (1-\varphi_0(V_0))\, V_0 \le (1-\psi_0(V_0))\, V_0 = (1-\psi_0(U_0))\, U_0 = U_1\, . \label{ir1}
\end{equation}
Suppose we managed to prove, that first $i-1$ inequalities takes place:
\begin{equation}
V_j \le U_j \, , \qquad j=1,\, ...,\, i-1\, , \label{im1}
\end{equation}
then the following estimation can be achieved:
\begin{equation}
V_i = (1-\varphi_{i-1}(V_0,...,V_{i-1}))\, V_{i-1} \le (1-\psi_{i-1}(V_0,...,V_{i-1}))\, V_{i-1} 
\le (1-\psi_{i-1}(U_0,...,U_{i-1}))\, U_{i-1} = U_i\, . \label{i}
\end{equation}
The first step in the estimation above is a consequence of (\ref{Vi}), the second step -- consequence of (\ref{psi}), the third step -- consequence of (\ref{Monoton}) and the the forth -- consequence of (\ref{Ui}). Properties of the majorant (\ref{Lmajor}) are now consequence of (\ref{ir1}) and (\ref{i}) by induction. $\Box$
\end{pf}

\begin{rmk}
As matrix $A$ dimension grows, calculations using formula (\ref{eps}) become more and more complicated due to the large number of combinations. Nevertheless, in case every column of matrix $A$ contains more than half of non-zero elements, there does exist simple and effective estimate from below:
\begin{equation}
\epsilon_p(\bm A) \ge \hat{\epsilon}_p(\bm A) \equiv \min_{\varkappa = \overline{1,\,\lfloor n/2\rfloor}}
\left\lbrace 
\min_{j\in\cj}\sum\limits_{i\in\hat{\ci}(j,\varkappa)}
\frac{a_{ij}}{\sum\limits_{k=1}^n a_{kj}} +
\min_{j\in\cj}\sum\limits_{i\in\hat{\ci}(j,n-\varkappa)}
\frac{a_{ij}}{\sum\limits_{k=1}^n a_{kj}}
\right\rbrace 
\, , \label{epsge}
\end{equation}
where $\hat{\ci}(j,\varkappa)$ denotes the set of indexes $i$, corresponding to the first $\varkappa$ minimal elements in the $j$ column of matrix $A$ (in particular, the following is true: $|\hat{\ci}(j,\varkappa)|=\varkappa$ for all $j\in\cj$, $\varkappa = \overline{1,\,n}$).
\end{rmk}

\begin{rmk} 
	As value $\epsilon_p(\bm A)$ (or $\hat{\epsilon}_p(\bm A)$) does not depend on the dimension of the problem, convergence speed of the method (\ref{s}-\ref{sp1}) might be high even for very large-dimensional problems. Convergence region (\ref{OblSx}) is constant due to the same reason. Stability of the method to the growth of $n$ and $m$ is one of its largest advantages.
	\end{rmk}

\begin{rmk} \label{2variants}
	In fact, two-step procedure (\ref{s}-\ref{sp1}) has two different realizations, not equivalent to each other. To enlarge the method's chances for global convergence,  (\ref{sp1}) should be chosen as a first step of the iteration instead of (\ref{s}), if for the right part of (\ref{OblSx}) the following inequality holds:
	
\begin{equation}
\epsilon_p(\bm A)\min_{i\in\ci} p_i < \epsilon_q(\bm A)\min_{j\in\cj} q_j\, ,
\label{pmq}
\end{equation}
where
\begin{equation}
\epsilon_q(\bm A) \equiv 
\min_{\tilde{\ci}\subset\ci} \min_{\tilde{\cj}\subset\cj}\left\{
\min_{i\in\tilde{\ci}}\sum\limits_{j\in\cj\setminus\tilde{\cj}}\frac{a_{ij}}{\sum\limits_{l=1}^m a_{il}} +
\min_{i\in\ci\setminus\tilde{\ci}}\sum\limits_{j\in\tilde{\cj}}\frac{a_{ij}}{\sum\limits_{l=1}^m a_{il}}\right\}\, .
\label{epsq}
\end{equation}
\end{rmk}

\begin{rmk} 
In case direct calculation using formula (\ref{eps}) (or (\ref{epsq})) can be performed or (\ref{epsge}) provides with not too low quality estimate, the a priori convergence criterion can be checked with higher precision using $r$-function (\ref{RF}). However, given the smallness of the relation under the "Big o"\ in the formula (\ref{OblSx}), the check of the criterion can be simplified to the last short formula.
\end{rmk}

\begin{rmk} Formula (\ref{Alpha}) provides us with the possibility of an a priori estimation of the number of iterations, necessary for the achievement of the desired accuracy. Thus we can easily estimate the total amount of calculations and expected calculation time. 
\end{rmk}

\begin{theorem}\ \ In case conditions of Theorem 1 are satisfied, iterations $\{\bm X^{(2t)}\}_{t=0}^{\infty}$ converge to the unique element of the admissible set, given by constraints (\ref{Xp}-\ref{Xq}).
\end{theorem}

\begin{pf} It's easy to check that if $\bm X\in\IR^{n\times m}$, then function
\begin{equation}
\|\bm X\|_1\equiv\sum_{j=1}^m \sum_{i=1}^n | x_{ij}|\, ,\label{Norma}
\end{equation}
is the norm in the respective vector space. Due to the method's recursive relationships we have:
\begin{eqnarray}
&&\|\bm X^{(s+1)} - \bm X^{(s)}\|_1 = \sum\limits_{j=1}^m\left|\sum\limits_{i=1}^n x^{(s)}_{ij} - q_j\,\right| = V(\bm X^{(s)})\, ,\\
&&\|\bm X^{(s)} - \bm X^{(s-1)}\|_1 = \sum\limits_{i=1}^n\left|\sum\limits_{j=1}^m x^{(s-1)}_{ij} - p_i\,\right| = V(\bm X^{(s-1)})\, .
\end{eqnarray}
As a consequence of these equalities, and taking into account the fact that Lyapunov function is not increasing, we get the following estimates:
\begin{equation}
\|\bm X^{(s+1)} - \bm X^{(s-1)}\|_1 \le 2\, V(\bm X^{(s-1)}) \le 
2 (1-\alpha)^{s-1} V(\bm X^{(0)})\, ,\qquad s=1, 3, 5, ...
\label{sosed}
\end{equation} 
The last step is a consequence of Lemma 1.

From (\ref{sosed}) follows, that all the iterations are situated inside a sphere of a finite radius (using the introduced norm) as a sum of decreasing geometric progression is finite. Thus we can conclude, that there does exist at least one accumulation point of iterations. However, if we suppose that there might be 
several different accumulation points $\bm X' \ne\bm X''$, it will contradict the fact that the norm of the difference between nearest iterations decreases with the speed of geometric progression, because in case of (\ref{sosed}) distance by norm between the limit points of the accumulation cannot be strictly more than zero. $\Box$
\end{pf}

\begin{cor}\ \ In case the conditions of Theorem 1 are satisfied and $V(\bm A) > 0$, the objective function on iterations $\{\bm X^{(2t)}\}_{t=0}^{\infty}$ is restricted from above by the number, which can be calculated a priori:
\begin{equation}
J(\bm X^{(2t)}) \le \frac{2\, V(\bm A)}{\alpha}\left(1-(1-\alpha)^{2t}\right)
< \frac{2\, V(\bm A)}{\alpha} \, ,\qquad t=0, 1, 2, ...
\label{Jup}
\end{equation}
where $\alpha$ is defined by formula (\ref{Alpha}) with
\begin{equation}
V_0 = V(\bm A)\, , \qquad \varepsilon_1 = \epsilon_p(\bm A)\quad\mbox{and}\quad
\varepsilon_2 = 1/\min_{i\in\ci} p_i\, .
\label{Foralpha}
\end{equation}
\end{cor}

\begin{pf} Since $\bm X^{(0)} = \bm A$,
\begin{eqnarray}
& \displaystyle J(\bm X^{(2t)})\equiv \|\bm X^{(2t)} - \bm A\|_2 \le \|\bm X^{(2t)} - \bm A\|_1
\le \sum\limits_{\tau = 1}^{t}\|\bm X^{(2\tau)} - \bm X^{(2\tau - 2)}\|_1 \le 
 2\, V(\bm A) \sum\limits_{\tau = 0}^{t-1} (1-\alpha)^{2\tau}\, . &
\end{eqnarray}
The last inequality is a consequence of (\ref{sosed}). $\Box$
\end{pf}

\begin{rmk} 
Algorithm (\ref{s}-\ref{sp1}) can be applied to QTP with the following heterogeneous constraints, which often arise in real-life tasks:
\begin{eqnarray}
& \displaystyle x_{ij} \ge d_{ij}\, ,\quad (i,j)\in\ci\times\cj\, ,\qquad
x_{ij} \le d_{ij}\, ,\quad (i,j)\in\cp\subset\ci\times\cj\, , &
\end{eqnarray}
under the following conditions:
\begin{eqnarray}
& \displaystyle p_i > \sum\limits_{j\in\cj} d_{ij} \, ,\quad i\in\ci\, , \qquad
q_j > \sum\limits_{i\in\ci} d_{ij}\, ,\quad j\in\cj\, , & \label{cond_pq}\\
& \displaystyle a_{ij} > d_{ij}\, ,\quad (i,j)\in(\ci\times\cj)\setminus\cp \, , \qquad
a_{ij} = d_{ij}\, ,\quad (i,j)\in\cp\subset\ci\times\cj\, .& \label{cond_a}
\end{eqnarray}
Indeed, after simple change of variables $x_{ij} = x'_{ij} + d_{ij}$ we have again the original task with homogeneous constraints, but with different initial data:
\begin{equation}
p'_i = p_i - \sum\limits_{j\in\cj} d_{ij}\, ,\quad
q'_j = q_j - \sum\limits_{i\in\ci} d_{ij}\, ,\quad
a'_{ij} = a_{ij} - d_{ij}\, .
\end{equation}
Conditions (\ref{cond_pq}-\ref{cond_a}) are nothing else than application of the basic assumptions of the task to the initial data after change of variables. If these conditions are satisfied, convergence criterion (\ref{OblSx}) will also work for the data, extended by matrix  $\left(d_{ij}\right)$.
\end{rmk}

\section{Estimation of the relative error by the objective function} 
\noindent
In case we omit in the original QTP all the restrictions, except from transport type restrictions (\ref{Xp}-\ref{Xq}) themselves, quadratic programming problem will be reduced to the system of linear algebraic equations which can be solved analytically. As a solution to this simplified problem we have (see Appendix):
\begin{equation}
\hat x_{ij} = a_{ij} - 
\frac{1}{|\cj|}\left(\sum\limits_{j\in\cj} a_{ij} - p_i\right) - \frac{1}{|\ci|}\left(\sum\limits_{i\in\ci} a_{ij} - q_j\right) + \frac{1}{|\ci||\cj|}\left(\sum\limits_{i\in\ci}\sum\limits_{j\in\cj} a_{ij} - s\right) . 
\label{ikshat}
\end{equation}
Now it is easy to receive the following estimates for the optimal value by the objective function:
\begin{equation}
J(\hat{\bm X}) \le J(\bm{X}^{*}) \le J(\bm X^{(s)})\, , 
\label{ZelOU}
\end{equation}
where $\bm X^{*}$ is the exact solution of the QTP with additional constraints (which is also unique). If $V(\bm A) > 0$, inequality (\ref{ZelOU}) provides us with the estimate of the relative error $\delta J$ of the solution, obtained by the method, compared with the exact value of the objective function in the optimal point $\bm X^{*}$:
\begin{equation}
\delta J(\bm X^{(s)})\equiv\frac{J(\bm X^{(s)})-J(\hat{\bm X})}{J(\hat{\bm X})} \ge
\frac{J(\bm X^{(s)})-J(\bm{X}^{*})}{J(\bm{X}^{*}\vphantom{\hat{\bm X}})} \, . 
\label{ZelOtnos}
\end{equation}

Unlike the error by restrictions $V(\bm X^{(s)})$, objective function error $\delta J(\bm X^{(s)})$ does not have to approach zero even in case when conditions of the Theorem 1 are satisfied. Nevertheless, there are two important points to mention here. First point is that relative error by the objective function can be a posteriori measured, which means that we have a "common sense filter". Second point is that when
the error by restrictions become neglectable after sufficient number of iterations, solution's error has to be small in case relative error by the objective function is small. In particular, the following asymptotic estimate is true: $\|\bm X^{(s)} - \bm X^{*}\|_2 \lesssim J(\bm X^{*})\sqrt{2\,\delta J(\bm X^{(s)})}$ if $\delta J(\bm X^{(s)}) \ll 1$.

\begin{cor}\ \ In case the conditions of Theorem 1 are satisfied and $V(\bm A) > 0$, the distance between iterations $\{\bm X^{(2t)}\}_{t=0}^{\infty}$ and  exact solution of problem (\ref{J}-\ref{Dop}) is restricted by the following expression:
\begin{equation}
\|\bm X^{(2t)} - \bm X^{*}\|_2 \le \sqrt{J^2(\bm X^{*})\delta J(\bm X^{(2t)})(2+
\delta J(\bm X^{(2t)})) + J(\bm X^{*})\,\frac{4\, V(\bm A)}{\alpha}(1-\alpha)^{2t}} \, ,\quad t=0, 1, 2, ...
\label{dXup}
\end{equation}
where $\alpha$ is defined by formula (\ref{Alpha}) with
\begin{equation}
V_0 = V(\bm A)\, , \qquad \varepsilon_1 = \epsilon_p(\bm A)\quad\mbox{and}\quad
\varepsilon_2 = 1/\min_{i\in\ci} p_i\, . \nonumber
\end{equation}
\end{cor}
\begin{pf} From  Theorem 2 follows, that there does exist a unique accumulation point of iterations $\bm X'$, and, besides, $V(\bm X') = 0$ due to the continuity of the chosen Lyapunov function. Thereby given limit point satisfy full system of restrictions of the problem (\ref{J}-\ref{Dop}), and for it the following inequality holds:
\begin{equation}
\left(\nabla J^2(\bm X^{*}),\,\bm X' - \bm X^{*} \right)_2 \ge 0\, ,
\label{Polayk1}
\end{equation}
which is true for the case of minimization of the convex differentiable function subjected to the convex restrictions \cite{Polyak}. As gradient can be calculated explicitly, inequality (\ref{Polayk1}) is equivalent to the following:
\begin{equation}
\left(\bm X^{*} - \bm A,\,\bm X' - \bm X^{*} \right)_2 \ge 0\, .
\label{Polayk2}
\end{equation}
Let's now write down the following evident equality, equivalent to the simple binomial expansion:
\begin{equation}
\|\bm X^{(2t)} - \bm A\|_2^2 = \|\bm X^{*} - \bm A\|_2^2 +
2 \left(\bm X^{*} - \bm A,\,\bm X^{(2t)} - \bm X^{*} \right)_2 +
\|\bm X^{(2t)} - \bm X^{*}\|_2^2 \, ,
\label{Tozhdestvo}
\end{equation}
we will now regroup the terms and make an estimate from above using (\ref{Polayk2}) and Cauchy inequality:
\begin{eqnarray}
& \displaystyle \|\bm X^{(2t)} - \bm X^{*}\|_2^2 = J^2(\bm X^{(2t)}) - J^2(\bm X^{*}) 
- 2 \left(\bm X^{*} - \bm A,\,\bm X' - \bm X^{*} \right)_2 
+ 2 \left(\bm X^{*} - \bm A,\,\bm X' - \bm X^{(2t)}\right)_2 \le {} & \nonumber\\
& \displaystyle {} \le J^2(\bm X^{(2t)}) - J^2(\bm X^{*}) + 2\,\|\bm X^{*} - \bm A\|_2\,\|\bm X' - \bm X^{(2t)}\|_2 \le {}
& \nonumber\\
& \displaystyle {} \le J^2(\bm X^{(2t)}) - J^2(\bm X^{*}) + 2\,J(\bm X^{*})\,\|\bm X' - \bm X^{(2t)}\|_1\le {}
& \nonumber\\
& \displaystyle {} \le J^2(\bm X^{*})\left((1 + \delta J(\bm X^{(2t)}))^2 - 1\right) + 2\,J(\bm X^{*})\,\sum\limits_{\tau = t}^{\infty}
\|\bm X^{(2\tau + 2)} - \bm X^{(2\tau )}\|_1 \, . & \nonumber
\end{eqnarray}
After some obvious simplifications and using (\ref{sosed}) for every member of the infinite sum we get (\ref{dXup}). $\Box$
\end{pf}
\pagebreak

\section{Numerical experiments}
\noindent
Both solution's properties and algorithm's behaviour depends on input data and parameters of the task. A number of experiments were conducted to study algorithm's behaviour dependence on the following factors:
\begin{itemize}
    \item total task dimension $n\times m$
    \item dimensions ratio $m/n$
    \item percent of zeros in matrix $\bm A$ (number of additional restrictions)
    \item difference in the values of elements of matrix $\bm A$ and vectors $\bm p$ and $\bm q$ (\textit{turbulence}).
\end{itemize}
Experiments were carried out for matrices of two types: one with ratio $m/n\simeq 1.7$, another with ratio $m/n\simeq 3.5$. Their results are presented in Table~1 and Table~2. To ensure results comparability, gene\-rator of random matrices and vectors created an input data in such a way, that initial (as well as the final) value of the Lyapunov function was the same in all experiments. This is also true about value of $s$ (sum of the elements of vectors). Moreover, random matrix generator can be adjusted to maintain certain average  percentage of zeros in matrix $\bm A$, as well as degree of turbulence in the elements values. For every experiment both variations of the method were performed according to Remark~\ref{2variants}.

The main result of experiments is confirmation of a hypothesis that the number of iterations necessary for convergence does not increase with the growth of the task's total dimension. We know that the method converges with the speed of geometric progression, but while dimension of the problem increases, ratio of this progression does not approach 1. This property of the method makes it possible to solve (approximately) QTPs with total dimension of more than a million using an ordinary home PC without parallelization, and an approximate solution time for such a problem is about half an hour.

Numerical experiments have proven that the method is stable to the growth of the dimensions ratio $m/n$ as well as to the growth of the turbulence in the initial data. When the percentage of zeros increases, the convergence of the method deteriorates. In the tables we compare the necessary number of iterations for the cases of $\sim\!\! 7$\% and $\sim\!\! 25$\% of zeros in matrix $\bm A$. In all the experiments only formula (\ref{epsge}) was used, so we had not tested the case with the percentage of zeros of about $50$\% or more, as the objective was to observe non-zero values of the lower estimates. On the other hand, the large number of zeros leads to a high probability of the system of restrictions to be inconsistent. That is why performing numerical experiments with high percentage of zeros in  matrix $\bm A$ is far from easy.

Another interesting result of the experiment is that relative error $\delta J$ is not too high, even for the largest dimensions studied. That explains why the method preserves the basic properties of the QTP solution, although the heuristic "does not know"\ which norm is chosen for objective function $J$. E.g., when average relative correction of all elements of the original matrix $\overline{(x-a)/a}$ is small, there are only a few elements that are corrected significantly, which is typical for the Euclidean norm.
\begin{landscape}
\begin{table}
    \centering
\footnotesize \begin{tabular}{|c||c|c||c|c||c|c|}
\hline
$n$ & $2$ & $2$ & $35$ & $35$ & $600$ & $600$ \\[0.pt]
\hline
$m$ & $7$ & $7$ & $121$ & $121$ & $2100$ & $2100$ \\[0.pt]
\hline
$n\times m$ & $14$ & $14$ & $4235$ & $4235$ & $1\,260\,000$ & $1\,260\,000$ \\[0.pt]
\hline
\hline
number of tests & $25$ & $25$ & $15$ & $15$ & $5$ & $5$ \\[0.pt]
\hline
$s_f$ & $27.66 \pm 9.47$ & $16.54 \pm 8.63$ & $7.16 \pm 0.63$ & $4.93 \pm 0.24$ & $7.00 \pm 0.00$ & $4.00 \pm 0.00$ \\[0.pt]
\hline
$\max s_f$ & $40$ & $35$ & $9$ & $5$ & $7$ & $4$ \\[0.pt]
\hline
$\sum a_{ij}$ & $(52.73 \pm 0.01)\,10^5$ & $(52.90 \pm 0.20)\,10^5$ & $(52.82 \pm 0.09)\,10^5$ & $(53.09 \pm 0.07)\,10^5$ & $(53.65 \pm 0.00)\,10^5$ & $(53.70 \pm 0.00)\,10^5$ \\[0.pt]
\hline
zeros (\%)& $24.85 \pm 3.56 $ & $ 6.85 \pm 1.39 $ & $25.01 \pm 0.02 $ & $ 6.98 \pm 0.43 $ & $24.99 \pm 0.16 $ & $ 7.11 \pm 0.07 $ \\[0.pt]
\hline
$a_{\max}/a_{\min}$ & $18.59 \pm 16.78$ & $28.81 \pm 24.64$ & $12189.37 \pm 7428.60$ & $12126.50 \pm 7130.41$ & $172.20 \pm 30.41$ & $178.00 \pm 28.63$ \\[0.pt]
\hline
$p_{\min}$ & $(1810.63 \pm 485.99)\,10^3$ & $(1698.46 \pm 547.18)\,10^3$ & $(3594.50 \pm 453.01)\,10$ & $(3414.27 \pm 390.49)\,10$ & $1789.40 \pm 17.47$ & $1788.80 \pm 17.90$ \\[0.pt]
\hline
$p_{\max}/p_{\min}$ & $2.27 \pm 1.36$ & $2.62 \pm 1.65$ & $13.93 \pm 2.78$ & $16.78 \pm 2.63$ & $28.59 \pm 1.96$ & $28.72 \pm 1.98$ \\[0.pt]
\hline
$q_{\min}$ & $(245.29 \pm 68.01)\,10^3$ & $(238.16 \pm 49.86)\,10^3$ & $9135.06 \pm 473.16$ & $9157.20 \pm 291.49$ & $508.20 \pm 0.40$ & $508.00 \pm 0.00$ \\[0.pt]
\hline
$q_{\max}/q_{\min}$ & $7.38 \pm 3.16$ & $8.23 \pm 3.38$ & $23.27 \pm 4.77$ & $22.03 \pm 4.30$ & $35.16 \pm 4.62$ & $36.01 \pm 4.13$ \\[0.pt]
\hline
$Z_p(X)$ & $(20.35 \pm 99.69)\,10^{-3}$ & $(378.47 \pm 217.36)\,10^{-3}$ & $(22.82 \pm 24.36)\,10^{-4}$ & $(240.39 \pm 41.25)\,10^{-3}$ & $(15.81 \pm 0.48)\,10^{-4}$ & $(295.69 \pm 7.54)\,10^{-3}$ \\[0.pt]
\hline
$\min Z_p(X)$ & ------ & $76.44\,10^{-3}$ & $ 0.82\,10^{-4}$ & $149.69\,10^{-3}$ & $15.33\,10^{-4}$ & $281.76\,10^{-3}$ \\[0.pt]
\hline
$Z_q(X)$ & $(197.15 \pm 97.29)\,10^{-3}$ & $(348.69 \pm 92.51)\,10^{-3}$ & $(138.39 \pm 173.82)\,10^{-4}$ & $(312.65 \pm 24.25)\,10^{-3}$ & $(-3.05 \pm 0.72)\,10^{-4}$ & $(320.59 \pm 16.15)\,10^{-3}$ \\[0.pt]
\hline
$\min Z_q(X)$ & $69.41\,10^{-3}$ & $207.47\,10^{-3}$ & $ 0.61\,10^{-4}$ & $272.90\,10^{-3}$ & $-3.93\,10^{-4}$ & $301.31\,10^{-3}$ \\[0.pt]
\hline
$\hat\epsilon_p(X)$ & $(52.37 \pm 124.55)\,10^{-3}$ & $(378.47 \pm 217.36)\,10^{-3}$ & $(23.11 \pm 24.31)\,10^{-4}$ & $(240.42 \pm 41.25)\,10^{-3}$ & $(18.14 \pm 0.46)\,10^{-4}$ & $(295.85 \pm 7.54)\,10^{-3}$ \\[0.pt]
\hline
$\max\hat\epsilon_p /\min\hat\epsilon_p$ & ------ & $11.14$ & $70.80$ & $2.19$ & $1.07$ & $1.07$ \\[0.pt]
\hline
$\hat\epsilon_q(X)$ & $(197.15 \pm 97.29)\,10^{-3}$ & $(348.69 \pm 92.51)\,10^{-3}$ & $(139.51 \pm 174.07)\,10^{-4}$ & $(312.70 \pm 24.23)\,10^{-3}$ & $(5.27 \pm 0.26)\,10^{-4}$ & $(321.27 \pm 16.06)\,10^{-3}$ \\[0.pt]
\hline
$\max\hat\epsilon_q /\min\hat\epsilon_q$ & $5.23$ & $2.72$ & $252.58$ & $1.30$ & $1.16$ & $1.12$ \\[0.pt]
\hline
$\delta J(X)$ (\%)& $55.72 \pm 15.42$ & $59.50 \pm 13.85$ & $30.71 \pm 05.79$ & $27.74 \pm 03.91$ & $48.54 \pm 01.08$ & $25.56 \pm 01.51$ \\[0.pt]
\hline
$\overline{(x-a)/a}$ (\%)& $1.93 \pm 1.12$ & $3.37 \pm 1.17$ & $1.52 \pm 0.47$ & $2.26 \pm 0.14$ & $3.21 \pm 0.07$ & $2.09 \pm 0.03$ \\[0.pt]
\hline
$\max (x-a)/a$ (\%)& $4.49 \pm 2.54$ & $10.94 \pm 6.20$ & $7.75 \pm 3.07$ & $12.97 \pm 2.31$ & $37.22 \pm 1.36$ & $11.90 \pm 0.36$ \\[0.pt]
\hline
calcul. time (s)& ${} < 10^{-2}$ & ${} < 10^{-2}$ & $0.33 \pm 0.05$ & $0.24 \pm 0.03$ & $5118.61 \pm 550.44$ & $2639.52 \pm 133.45$ \\[0.pt]
\hline
\end{tabular}
\vspace{5.pt}
    \caption{Dependence of the algorithm's behavior and the solution's properties on the problem's dimension for two sets of additional constraints (percent of zeros). Case $m/n \simeq 3.5$ (far from square matrix), $V(X^{(0)})\simeq 10^5$, $V(X^{\!f})\simeq 1$ and $s = 5320851$.}
    \label{TABLE4}
\end{table}
\begin{table}
    \centering
\footnotesize \begin{tabular}{|c||c|c||c|c||c|c|}
\hline
$n$ & $3$ & $3$ & $52$ & $52$ & $900$ & $900$ \\[0.pt]
\hline
$m$ & $5$ & $5$ & $87$ & $87$ & $1500$ & $1500$ \\[0.pt]
\hline
$n\times m$ & $15$ & $15$ & $4524$ & $4524$ & $1\,350\,000$ & $1\,350\,000$ \\[0.pt]
\hline
\hline
number of tests & $25$ & $25$ & $15$ & $15$ & $5$ & $5$ \\[0.pt]
\hline
$s_f$ & $25.06 \pm 16.94$ & $9.66 \pm 3.40$ & $7.13 \pm 0.66$ & $5.00 \pm 0.00$ & $7.00 \pm 0.00$ & $4.00 \pm 0.00$ \\[0.pt]
\hline
$\max s_f$ & $79$ & $20$ & $8$ & $5$ & $7$ & $4$ \\[0.pt]
\hline
$\sum a_{ij}$ & $(55.95 \pm 0.18)\,10^5$ & $(56.13 \pm 0.19)\,10^5$ & $(56.02 \pm 0.10)\,10^5$ & $(56.23 \pm 0.06)\,10^5$ & $(56.80 \pm 0.00)\,10^5$ & $(56.85 \pm 0.00)\,10^5$ \\[0.pt]
\hline
zeros (\%)& $24.80 \pm 2.99 $ & $ 7.20 \pm 1.80 $ & $25.09 \pm 0.18 $ & $ 7.11 \pm 0.35 $ & $24.93 \pm 0.23 $ & $ 7.08 \pm 0.16 $ \\[0.pt]
\hline
$a_{\max}/a_{\min}$ & $29.35 \pm 25.95$ & $20.47 \pm 8.89$ & $11922.16 \pm 5062.84$ & $11281.26 \pm 6718.33$ & $159.00 \pm 19.04$ & $164.60 \pm 62.79$ \\[0.pt]
\hline
$p_{\min}$ & $(751.65 \pm 303.48)\,10^3$ & $(886.04 \pm 329.34)\,10^3$ & $(2442.27 \pm 230.52)\,10$ & $(2338.86 \pm 185.42)\,10$ & $1260.40 \pm 7.49$ & $1257.60 \pm 1.01$ \\[0.pt]
\hline
$p_{\max}/p_{\min}$ & $5.32 \pm 2.80$ & $4.09 \pm 2.29$ & $14.95 \pm 2.70$ & $18.25 \pm 5.58$ & $30.90 \pm 2.07$ & $29.46 \pm 6.02$ \\[0.pt]
\hline
$q_{\min}$ & $(403.67 \pm 134.73)\,10^3$ & $(384.92 \pm 114.28)\,10^3$ & $(132.46 \pm 1.97)\,10^2$ & $(134.05 \pm 4.22)\,10^2$ & $753.00 \pm 1.26$ & $752.60 \pm 0.48$ \\[0.pt]
\hline
$q_{\max}/q_{\min}$ & $6.02 \pm 2.68$ & $6.50 \pm 2.39$ & $18.86 \pm 3.17$ & $21.23 \pm 4.82$ & $30.52 \pm 2.15$ & $32.34 \pm 4.98$ \\[0.pt]
\hline
$Z_p(X)$ & $(188.76 \pm 89.36)\,10^{-3}$ & $(310.61 \pm 119.44)\,10^{-3}$ & $(7.64 \pm 6.70)\,10^{-4}$ & $(265.35 \pm 34.85)\,10^{-3}$ & $(9.44 \pm 0.24)\,10^{-4}$ & $(307.08 \pm 13.26)\,10^{-3}$ \\[0.pt]
\hline
$\min Z_p(X)$ & $28.25\,10^{-3}$ & $85.83\,10^{-3}$ & $ 0.56\,10^{-4}$ & $221.62\,10^{-3}$ & $9.10\,10^{-4}$ & $284.10\,10^{-3}$ \\[0.pt]
\hline
$Z_q(X)$ & $(154.56 \pm 65.82)\,10^{-3}$ & $(360.99 \pm 129.24)\,10^{-3}$ & $(8.85 \pm 14.03)\,10^{-4}$ & $(288.97 \pm 23.09)\,10^{-3}$ & $(2.83 \pm 0.27)\,10^{-4}$ & $(313.39 \pm 16.25)\,10^{-3}$ \\[0.pt]
\hline
$\min Z_q(X)$ & $41.22\,10^{-3}$ & $148.01\,10^{-3}$ & $-0.48\,10^{-4}$ & $230.39\,10^{-3}$ & $2.43\,10^{-4}$ & $290.03\,10^{-3}$ \\[0.pt]
\hline
$\hat\epsilon_p(X)$ & $(188.77 \pm 89.36)\,10^{-3}$ & $(310.61 \pm 119.44)\,10^{-3}$ & $(8.04 \pm 6.71)\,10^{-4}$ & $(265.37 \pm 34.85)\,10^{-3}$ & $(12.07 \pm 0.13)\,10^{-4}$ & $(307.21 \pm 13.26)\,10^{-3}$ \\[0.pt]
\hline
$\max\hat\epsilon_p /\min\hat\epsilon_p$ & $12.84$ & $5.76$ & $24.77$ & $1.67$ & $1.03$ & $1.12$ \\[0.pt]
\hline
$\hat\epsilon_q(X)$ & $(154.57 \pm 65.81)\,10^{-3}$ & $(360.99 \pm 129.24)\,10^{-3}$ & $(9.72 \pm 13.91)\,10^{-4}$ & $(288.99 \pm 23.09)\,10^{-3}$ & $(7.53 \pm 0.19)\,10^{-4}$ & $(313.58 \pm 16.16)\,10^{-3}$ \\[0.pt]
\hline
$\max\hat\epsilon_q /\min\hat\epsilon_q$ & $8.02$ & $3.93$ & $164.33$ & $1.39$ & $1.07$ & $1.16$ \\[0.pt]
\hline
$\delta J(X)$ (\%)& $55.46 \pm 15.07$ & $53.00 \pm 14.12$ & $31.94 \pm 7.01$ & $29.25 \pm 7.22$ & $54.04 \pm 0.54$ & $23.71 \pm 0.61$ \\[0.pt]
\hline
$\overline{(x-a)/a}$ (\%)& $2.04 \pm 1.05$ & $2.67 \pm 0.48$ & $1.68 \pm 0.39$ & $2.13 \pm 0.27$ & $3.27 \pm 0.02$ & $1.88 \pm .02$ \\[0.pt]
\hline
$\max (x-a)/a$ (\%)& $4.15 \pm 2.15$ & $7.10 \pm 2.96$ & $13.11 \pm 5.78$ & $14.94 \pm 2.55$ & $37.67 \pm 3.91$ & $11.67 \pm 0.40$ \\[0.pt]
\hline
calcul. time (s)& ${} < 10^{-2}$ & ${} < 10^{-2}$ & $0.35 \pm 0.04$ & $0.27 \pm 0.02$ & $1940.98 \pm 354.74$ & $3371.94 \pm 1027.59$ \\[0.pt]
\hline
\end{tabular}
\vspace{5.pt}
    \caption{Dependence of the algorithm's behavior and the solution's properties on the problem's dimension
    for two sets of additional constraints (percent of zeros). Case $m/n \simeq 1.7$ (close to square matrix),
   $V(X^{(0)})\simeq 10^5$, $V(X^{\!f})\simeq 1$ and $s = 5635083$.}
    \label{TABLE5}
\end{table}
\begin{table}
    \centering
\vspace{10.pt}
\normalsize \begin{tabular}{|c||c|c|c||c|c|}
\hline
$n$ & $52$ & $35$ & $52$ & $52$ & $52$ \\[0.pt]
\hline
$m$ & $87$ & $121$ & $87$ & $87$ & $87$ \\[0.pt]
\hline
$n\times m$ & $4524$ & $4235$ & $4524$ & $4524$ & $4524$ \\[0.pt]
\hline
\hline
number of tests & $15$ & $15$ & $15$ & $15$ & $15$ \\[0.pt]
\hline
$s_f$ & $6.80 \pm 0.65$ & $6.86 \pm 0.56$ & $4.06 \pm 0.24$ & $7.13 \pm 0.42$ & $6.66 \pm 0.59$ \\[0.pt]
\hline
$\max s_f$ & $8$ & $8$ & $5$ & $8$ & $8$ \\[0.pt]
\hline
$\sum a_{ij}$ & $(56.02 \pm 0.10)\,10^5$ & $(52.82 \pm 0.09)\,10^5$ & $(56.23 \pm 0.06)\,10^5$ & $(55.85 \pm 0.00)\,10^5$ & $(55.85 \pm 0.00)\,10^5$ \\[0.pt]
\hline
zeros (\%)& $25.09 \pm 0.18 $ & $25.01 \pm 0.02 $ & $ 7.11 \pm 0.35 $ & $25.00 \pm 0.00 $ & $25.00 \pm 0.00 $ \\[0.pt]
\hline
$a_{\max}/a_{\min}$ & $11922.16 \pm 5062.84$ & $12189.37 \pm 7428.60$ & $11281.26 \pm 6718.33$ & $1278.50 \pm 937.13$ & $1278.50 \pm 937.13$ \\[0.pt]
\hline
$p_{\min}$ & $(2442.27 \pm 230.51)\,10$ & $(3594.50 \pm 453.01)\,10$ & $(2338.86 \pm 185.42)\,10$ & $(6585.46 \pm 73.73)\,10$ & $(6585.46 \pm 73.73)\,10$ \\[0.pt]
\hline
$p_{\max}/p_{\min}$ & $14.95 \pm 2.70$ & $13.93 \pm 2.78$ & $18.25 \pm 5.58$ & $4.01 \pm 0.81$ & $4.01 \pm 0.81$ \\[0.pt]
\hline
$q_{\min}$ & $(1324.64 \pm 19.69)\,10$ & $(913.51 \pm 47.32)\,10$ & $(1340.49 \pm 42.24)\,10$ & $(3906.12 \pm 21.32)\,10$ & $(3906.12 \pm 21.32)\,10$ \\[0.pt]
\hline
$q_{\max}/q_{\min}$ & $18.86 \pm 3.17$ & $23.27 \pm 4.77$ & $21.23 \pm 4.82$ & $4.16 \pm 0.76$ & $4.16 \pm 0.76$ \\[0.pt]
\hline
$Z_p(X)$ & $( 0.22 \pm 2.30)\,10^{-4}$ & $(19.17 \pm 26.19)\,10^{-4}$ & $(265.71 \pm 35.64)\,10^{-3}$ & $(82.62 \pm 30.67)\,10^{-3}$ & $(84.54 \pm 30.86)\,10^{-3}$ \\[0.pt]
\hline
$\min Z_p(X)$ & $-0.84\,10^{-4}$ & $-0.26\,10^{-4}$ & $217.18\,10^{-3}$ & $36.75\,10^{-3}$ & $38.57\,10^{-3}$ \\[0.pt]
\hline
$\hat\epsilon_p(X)$ & $( 0.60 \pm 2.26)\,10^{-4}$ & $(19.41 \pm 26.23)\,10^{-4}$ & $(265.75 \pm 35.63)\,10^{-3}$ & $(82.63 \pm 30.67)\,10^{-3}$ & $(84.56 \pm 30.86)\,10^{-3}$ \\[0.pt]
\hline
new zeros & $  7.53 \pm 3.49$ & $  7.20 \pm 5.85$ & $  6.00 \pm 3.82$ & $0$ & $0$ \\[0.pt]
\hline
$J(X)$& $1648.22 \pm 150.63$ & $1666.28 \pm 153.39$ & $1468.98 \pm 37.62$ & $1239.76 \pm 36.91$ & $1175.50 \pm 30.15$ \\[0.pt]
\hline
$\delta J(X)$ (\%)& $15.72 \pm 3.80$ & $14.78 \pm 1.90$ & $ 2.54 \pm 0.40$ & $23.61 \pm 0.94$ & $17.22 \pm 0.66$ \\[0.pt]
\hline
$\overline{(x-a)/a}$ (\%)& $3.07 \pm 0.81$ & $2.78 \pm 1.02$ & $3.73 \pm 0.64$ & $ 0.90 \pm 0.02$ & $1.13 \pm 0.05$ \\[0.pt]
\hline
$\max (x-a)/a$ (\%)& $370.88 \pm 245.10$ & $309.24 \pm 379.52$ & $322.07 \pm 284.88$ & $3.34 \pm 0.27$ & $109.55 \pm 116.57$ \\[0.pt]
\hline
calcul. time (s)& $0.45 \pm 0.66$ & $0.31 \pm 0.03$ & $0.21 \pm 0.04$ & $0.36 \pm 0.06$ & $0.42 \pm 0.09$ \\[0.pt]
\hline
\end{tabular}
\vspace{10.pt}
    \caption{Dependence of the algorithm's behavior and the solution's properties on the use of combined approach under various problem's parameters, $V(X^{(0)})\simeq 10^5$ and $V(X^{\!f})\simeq 1$.}
    \label{TABLE6}
\end{table}
\end{landscape}

The result quality can be improved using combined approach, which consists of the following steps:

\begin{itemize}
    \item  using formula (\ref{ikshat}), exact analytical solution $\hat{\bm X}$ of the problem with less constraints is calculated 
    \item to satisfy restrictions (\ref{Neotric}-\ref{Dop}), all negative elements as well as all elements in positions  $(i,j)\in\cp$ are replaced by zeros
    \item new matrix, obtained after described calculations, is used instead of matrix $\bm A$ as a starting point in the iterative procedure.
\end{itemize}

To study the method's reaction to this combined approach, additional experiments were performed (see Table~3). Results of the experiments using combined approach can be found in columns 1, 2, 3 and 5 of the table. It is noteworthy that when turbulence is high (columns 1-3) new zeros appear in the solution matrix. Most significant effect of the combined approach is the fact that estimation of the relative error $\delta J$ at the final iteration is reduced by approximately a factor of two.

The method's convergence investigation was performed by monitoring $\hat\epsilon_p(\bm X^{(s)})$,  $\hat\epsilon_q(\bm X^{(s)})$ and residuals:
\begin{equation}
    Z_p(\bm X^{(s)})\equiv \hat\epsilon_p(\bm X^{(s)}) -
    \frac{2\,V(\bm X^{(s)})}{\min_{i\in\ci} p_i}
		\left(1+O\!\left(\frac{V(\bm X^{(s)})}{\min_{i\in\ci} p_i}\right)\right)\,  
\end{equation}
and $Z_q(\bm X^{(s)})$ (formula is similar), calculated based on formulas (\ref{OblSx}) and (\ref{Zona2}) with $\varepsilon_1 = \hat\epsilon_p$ or $\varepsilon_1 = \hat\epsilon_q$ and $\varepsilon_2 = 1/\min_{i\in\ci} p_i$ or $\varepsilon_2 = 1/\min_{j\in\cj} q_j$. Positive values of $Z_p$ or $Z_q$ a posteriori prove that the full system of restrictions (\ref{Xp}-\ref{Dop}) is consistent. The method's convergence is rapidly deteriorating when restrictions are close to inconsistent.

\section{Conclusions}

\noindent
Quadratic transportation problem has been studied for quite a long time \cite{QTP94}. In spite of that, this task remains relevant \cite{QTP13} and has various practical applications \cite{QTP03}. Efforts were made to generalize QTP for the case of convex objective function with additional constraints \cite{QTP15}. In all these papers authors are trying to find exact solution of the transportation problem. Approach used in this article is to meet all the constraints precisely, and at the same time find an approximate solution by the objective function. This is especially valuable for forecasting purposes, as forecasts are never accurate due to the nature of forecasting itself. The fact that we are not obliged to find an exact solution allow the application of extremely effective computational algorithm,  thus letting us to solve transportation problems of very high dimension.

The problem of high dimension in the mathematical programming has been a challenge for quite a while. The main idea was to decompose the original problem. Key methods to be mentioned here are Dantzig-Wolfe decomposition \cite{DW1961}, Benders decomposition \cite{Bender62} and cross decomposition of Holmberg \cite{KH1994, KH1984}. Unlike the algorithm described in this paper, all these methods do not allow massive parallelization, solving problem of high dimension by longer calculation time. But even on a simple home PC without a hint of parallelization, the proposed method solves relatively large task in the acceptable time, which can be clearly seen from carried out numerical experiments.

Another thing we can see from the experiments, is that condition (\ref{OblSx}) for the convergence region is sufficient, but not necessary. In many practical cases, in spite of the fact that matrix $\bm X^{(0)}= \bm A$ does not satisfy condition (\ref{OblSx}), the algorithm described above does converge. Moreover, starting from a certain iteration, condition (\ref{OblSx}) is satisfied by the matrix of the current iteration, as well as by all the following ones. In case this does not happen and $ \epsilon(\bm X^{(s)}) = 0$ beginning from some iteration, the original matrix might allow block decomposition, and thus original problem can be decomposed into several problems of a smaller dimension with independent data. For every such a sub-problem the algorithm will converge if the corresponding sub-vectors $\bm p$ and $\bm q$ satisfy the balance condition, and in general, the corresponding complete sub-systems of constraints are compatible.

As we suppose that the original matrix $\bm A$ is pretty close to the admissible one, the number of active inequality constraints (\ref{Neotric}-\ref{Dop}) is relatively small. As a consequence, for the real-life cases formula (\ref{ZelOtnos}) provides us with quite an accurate estimate of the relative error by the objective function. Applying upper restriction for the value of $\delta J$, we can construct an a posteriori solution filter, and thus avoid any inappropriate solutions.

The method proposed most likely can be generalized to the multidimensional hierarchical time series, where the number of dimensions is more than 2. In that case matrix $\bm A$ will turn into a multi-index tensor. Convergence analysis of the algorithm will be more difficult, but similar by the structure to the one proposed above. The result on the convergence rate of the methods of the considered type, obtained in the article \cite{Goreinov}, let us hope that it will be possible to prove the linear convergence in more complicated cases. On the other hand, the same classical result explains why we cannot hope for a higher convergence rate even under additional assumptions.

\section{Appendix}

\noindent
Here we demonstrate output of formula (\ref{ikshat}). QTP without inequality constraints (\ref{Neotric}-\ref{Dop}) can be solved using the method of Lagrange multipliers. Applying this approach we get the following system of linear equations for primal ($x_{ij}$) and  dual ($\lambda_i$, $\mu_j$) variables:
\begin{eqnarray}
&& \displaystyle 
\sum\limits_{j=1}^m x_{ij} = p_i\, ,\quad i=\overline{1,\, n}\, , 
\nonumber\\
&& \displaystyle
\sum\limits_{i=1}^n x_{ij} = q_j\, ,\quad  j=\overline{1,\, m}\, , 
\nonumber\\[7.pt]
&&\displaystyle
x_{ij} = a_{ij} - (\lambda_i + \mu_j)\, ,\qquad i=\overline{1,\, n}\, ,\; j=\overline{1,\, m}
\label{dLdl}
\end{eqnarray}
Let's substitute the last equation into first two and exclude primal variables:
\begin{eqnarray}
&&\xi_i - m\,\lambda_i - \sum\limits_{j=1}^m \mu_j = p_i\, ,\qquad i=\overline{1,\, n}\, , \label{dual1}\\
&&\eta_j - n\,\mu_j - \sum\limits_{i=1}^n \lambda_i = q_j\, ,\qquad j=\overline{1,\, m}\, .\label{dual2}
\end{eqnarray}
where $\xi_i = \sum_{j=1}^m a_{ij}$ and $\eta_j = \sum_{i=1}^n a_{ij}$. Summing (\ref{dual1}) or (\ref{dual2}) by $i$ or by $j$ respectively, we obtain:
\begin{equation}
m \sum\limits_{i=1}^n \lambda_i + n \sum\limits_{j=1}^m \mu_j = \zeta - s
\label{decisive}
\end{equation}
with $\zeta = \sum_{j=1}^m\sum_{i=1}^n a_{ij}\,$, $s = \sum_{i=1}^n p_i = \sum_{j=1}^m q_j$. Let's now write down a system, equivalent to (\ref{dual1}-\ref{dual2}):
\begin{eqnarray}
&& \lambda_i = \frac{1}{m} \left(\xi_i - p_i - \sum\limits_{j=1}^m \mu_j\right)\, ,\qquad i=\overline{1,\, n}\, , \\
&& \mu_j = \frac{1}{n} \left(\eta_j - q_j - \sum\limits_{i=1}^n \lambda_i\right)\, ,\qquad j=\overline{1,\, m}\, .
\end{eqnarray}
Now we have:
\begin{equation}
\lambda_i + \mu_j =  \frac{1}{n\, m}\left( n\,(\xi_i - p_i) + m\,(\eta_j - q_j) - \left(m \sum\limits_{i=1}^n \lambda_i + n \sum\limits_{j=1}^m \mu_j\right)\right)\, ,\quad i=\overline{1,\, n}\, ,\; j=\overline{1,\, m}\, .
\label{LpM}
\end{equation}
If we substitute (\ref{decisive}) into (\ref{LpM}), we obtain a solution for the particular sum of the dual variables:
\begin{equation}
\lambda_i + \mu_j =  \frac{1}{n\, m}\left( n\,(\xi_i - p_i) + m\,(\eta_j - q_j) - (\zeta - s)\vphantom{2^2}\right)\, ,\qquad i=\overline{1,\, n}\, ,\: j=\overline{1,\, m}\, .
\label{dual_solution}
\end{equation}
After substitution of (\ref{dual_solution}) into (\ref{dLdl}), we finally receive desired analytical solution (\ref{ikshat}).

\end{document}